\documentclass[12pt]{amsart}
\usepackage{amssymb,latexsym,comment,url}

\newcommand{\EE}{{\mathcal E}}

\newcommand{\CC}{\operatorname{\mathcal C}}

\newcommand{\Q}{{\mathbb Q}}

\newcommand{\rank}{\operatorname{rank}}

\newcommand{\Z}{{\mathbb Z}}

\newcommand{\Magma}{{\sf MAGMA }}

\newenvironment{Proof}{\par\noindent{\sc Proof:}}%
                      {\hspace*{\fill}\nobreak$\Box$\par\medskip}
                       {\hspace*{\fill}\nobreak$\Box$\par\medskip}

\newtheorem{Proposition}{Proposition}[section]
\newtheorem{Theorem}[Proposition]{Theorem}

\newtheorem{Corollary}[Proposition]{Corollary}

\theoremstyle{definition}
\newtheorem{Definition}[Proposition]{Definition}
\newtheorem{Remark}[Proposition]{Remark}

\renewcommand{\baselinestretch}{1.1}

\begin{document}

\title[Sequences of consecutive squares on elliptic curves]%
{On sequences of consecutive squares on elliptic curves}

\author[M. Kamel]%
{Mohamed~Kamel}
\address{Department of Mathematics, Faculty of Science, Cairo University, Giza, Egypt}
\email{mohgamal@sci.cu.edu.eg}

\author[M. Sadek]%
{Mohammad~Sadek}
\address{American University in Cairo, Mathematics and Actuarial Science Department, AUC Avenue, New Cairo, Egypt}
\email{mmsadek@aucegypt.edu}

\date{}

\begin{abstract}
Let $C$ be an elliptic curve defined over $\Q$ by the equation $y^2=x^3+Ax+B$ where $A,B\in\Q$. A sequence of rational points $(x_i,y_i)\in C(\Q),\,i=1,2,\ldots,$ is said to form a sequence of consecutive squares on $C$ if the sequence of $x$-coordinates, $x_i,i=1,2,\ldots$, consists of consecutive squares. We produce an infinite family of elliptic curves $C$ with a $5$-term sequence of consecutive squares. Furthermore, this sequence consists of five independent rational points in $C(\Q)$. In particular, the rank $r$ of $C(\Q)$ satisfies $r\ge 5$.

\end{abstract}

\maketitle


\section{Introduction}
In \cite{Bremner}, Bremner initiated the discussion of certain arithmetic questions on rational points of elliptic curves attempting to relate the group structure on an elliptic curve $E$ to the addition group operation on the rational line. He raised the question of the existence of a sequence of rational points in $E(\Q)$ whose $x$-coordinates form an arithmetic progression in $\Q$. Such sequence is called an arithmetic progression sequence in $E(\Q)$. A variety of questions may be posed. For instance, how long these sequences can be and how many elliptic curves would have such long sequences of rational points. The existence of infinitely many elliptic curves with length 8 arithmetic progressions was proved.
Several authors introduced different approaches to find infinitely many elliptic curves with longer arithmetic progression sequences, see \cite{Alvarado, Campbell, Macleod, Ulas1}.

In \cite{BremnerUlas}, the study of sequences of rational points on elliptic curves whose $x$-coordinates form a geometric progression in $\Q$ was initiated. An infinite family of elliptic curves having geometric progression sequences of length 4 was exhibited. It was remarked that infinitely many elliptic curves with $5$-term geometric progression sequences can be constructed.

In this note, we discuss sequences of rational points on elliptic curves whose $x$-coordinates form a sequence of consecutive squares. We consider elliptic curves defined by the equation $y^2=ax^3+bx+c$ over $\Q$. We show that elliptic curves defined by the latter equation with 5-term sequences of rational points whose $x$-coordinates are elements in a sequence of consecutive squares in $\Q$ are parametrized by an elliptic surface whose rank is positive. Hence, one deduces the existence of infinitely many such elliptic curves. Moreover, we show that the five rational points forming the sequence are linearly independent in the group of rational points of the elliptic curve they lie on. In particular, we introduce an infinite family of elliptic curves of rank $\ge 5$.

\section{Sequences of Consecutive Squares}
\begin{Definition}
Let $C$ be an elliptic curve defined over a number field $K$ by the Weierstrass equation $y^2+a_1xy+a_3y=x^3+a_2x^2+a_4x+a_6,\,a_i\in K$. The sequence $(x_i,y_i)\in C(K)$ is said to be a {\em sequence of consecutive squares} on $C$ if there is a $u\in K$ such that $x_i=(u+i)^2$, $i=1,2,\ldots$.
\end{Definition}
The following proposition ensures the finiteness of the sequence of consecutive squares on an elliptic curve.
\begin{Proposition}
\label{prop1}
Let $C$ be an elliptic curve defined over a number field $K$ by a Weierstrass equation of the form
\[y^2+a_1xy+a_3y=x^3+a_2x^2+a_4x+a_6,\,a_i\in K.\]
Let $(x_i,y_i)\in C(K)$ be a sequence of consecutive squares on $C$. Then the sequence $(x_i,y_i)$ is finite.
\end{Proposition}
\begin{Proof}
We can assume without loss of generality that $x_i=(u+i)^2$, $i=1,2,\ldots$, $u\in K$. This sequence gives rise to a sequence of rational points on the genus $2$ hyperelliptic curve
\[\CC:y^2+a_1x^2y+a_3y=x^6+a_2x^4+a_4x^2+a_6.\]
Namely, the points $(u+i,y)\in \CC(K)$. According to Faltings' Theorem, \cite{Falting}, one knows that $\CC(K)$ is finite, hence the sequence is finite.
\end{Proof}
Based on the above proposition, one may present the following definition.
\begin{Definition}
Let $C$ be an elliptic curve over $\Q$ defined by a Weierstrass equation. Let $(x_i,y_i)\in C(\Q),\,i=1,2,\ldots,n,$ be a sequence of consecutive squares on $C$. Then $n$ is said to be the {\em length} of the sequence.
\end{Definition}
\section{Constructing elliptic curves with long sequences of consecutive squares}

In this note, we focus our attention on the family of elliptic curves given by the affine equation $C: y^2=a x^3+b x +c$ over $\Q$. We will show that there are infinitely many elliptic curves defined by the latter equation containing 5-term sequences of consecutive squares.

One observes that if $(t^2,d),((t+1)^2,e)$, and $((t+2)^2,f)$ lie in $C(\Q)$, where $t\in\Q$, then these rational points form a 3-term sequence of consecutive squares. Indeed, one has
\begin{eqnarray*}
d^2&=&at^6+bt^2+c\\e^2&=&a(t+1)^6+b(t+1)^2+c\\f^2&=&a(t+2)^6+b(t+2)^2+c.
\end{eqnarray*}
It is a standard linear algebra exercise to show that
{\footnotesize \begin{eqnarray}
\label{eqabc}
a&=&\frac{ (3 + 2t)d^2-4(1 + t)e^2+(1 + 2 t)f^2}{4 (15 + 73 t + 135 t^2 + 125 t^3 + 60 t^4 + 12 t^5)}\nonumber\\
b&=&\frac{(3+2 t) (3+3 t+t^2) (7+9 t+3 t^2)d^2+(1+t)(-4 (4+2 t+t^2) (4+6 t+3 t^2)e^2) }{ 4(1 + 2 t) (15 + 43 t + 49 t^2 + 27 t^3 + 6 t^4)}\nonumber\\
&+&\frac{ (1+t)(4+2 t+t^2) (4+6 t+3 t^2)f^2}{ 4(1 + 2 t) (15 + 43 t + 49 t^2 + 27 t^3 + 6 t^4)}\nonumber\\
c&=&\frac{(2 + t)^2 (15 + 43 t + 46 t^2 + 22 t^3 + 4 t^4)d^2  -8t^2(2 + t)^2 (2 + 2 t + t^2)e^2 +t^2 (1 + 5 t + 10 t^2 + 10 t^3 + 4 t^4)f^2}{4 (1 + 2 t) (15 + 28 t + 21 t^2 + 6 t^3)}.\nonumber\\
\end{eqnarray}}
In particular, one has the following result.
\begin{Remark}
\label{Rem1}
The above argument indicates that given $d,e,f\in\Q(t)$, there exist $a,b,c\in\Q(t)$ such that the ordered pairs $(t^2,d),((t+1)^2,e)$ and $((t+2)^2,f)$ are three rational points on the elliptic surface $y^2=ax^3+bx+c$.
\end{Remark}

Now, if $((t+3)^2,g)\in C(\Q)$, then one has a 4-term sequence of consecutive squares on $C$. In fact, using the above values for $a,b,c$, one then sees that
{\footnotesize\begin{eqnarray}
\label{eq1}
g^2=\frac{ (5 + 2 t)((2 + t) (14 + 12 t + 3 t^2)d^2 -3(1 + t) (13 + 10 t + 3 t^2)e^2)+3  (2 + t) (1 + 2 t) (10 + 8 t + 3 t^2)f^2 }{(1 + t) (1 + 2 t) (5 + 6 t + 3 t^2)}.\nonumber\\
\end{eqnarray}}
Therefore, in view of Remark \ref{Rem1}, one needs to find the elements $d,e,f$ and $g$ in $\Q(t)$ satisfying the latter equation in order to construct an elliptic curve $C$ with a 4-term sequence of consecutive squares. In fact, since $(d,e,f,g)=(1,1,1,1)$ is a solution for equation (\ref{eq1}), the general solution $(d,e,f,g)$ is given by the following parametrization:

{\footnotesize
\begin{eqnarray}\label{eq2}d&=& (2 + t) (5 + 2 t) (14 + 12 t + 3 t^2)p^2+3 (1 + t) (5 + 2 t) (13 + 10 t + 3 t^2)q^2-3 (2 + t) (1 + 2 t) (10 + 8 t + 3 t^2) w^2 \nonumber\\
      &-&6  (65+141 t+111 t^2+41 t^3+6 t^4)pq+6 (20+66 t+66 t^2+31 t^3+6 t^4)pw,\nonumber\\
     e&=& -(2 + t) (5 + 2 t) (14 + 12 t + 3 t^2)p^2-3  (1 + t) (5 + 2 t) (13 + 10 t + 3 t^2)q^2-3 (2 + t) (1 + 2 t) (10 + 8 t + 3 t^2) w^2 \nonumber\\
      &+&2  (140+246 t+166 t^2+51 t^3+6 t^4)p q+6  (20+66 t+66 t^2+31 t^3+6 t^4)q w,\nonumber\\
     f&=& -(2 + t) (5 + 2 t) (14 + 12 t + 3 t^2)p^2+3  (1 + t) (5 + 2 t) (13 + 10 t + 3 t^2)q^2+3 (2 + t) (1 + 2 t) (10 + 8 t + 3 t^2) w^2\nonumber\\
      &-&6 (1 + t) (5 + 2 t) (13 + 10 t + 3 t^2) q w+2  (2 + t) (5 + 2 t) (14 + 12 t + 3 t^2)p w,\nonumber\\
      g&=&-p^2(140+ 246t+ 166t^2+ 51t^3+6t^4)+3(q^2(65+ 141t+ 111t^2+ 41t^3+6t^4)\nonumber\\&-&(20+ 66t+ 66t^2+ 31t^3+6t^4)w^2).\nonumber\\
      \end{eqnarray}}
Consult \cite[\S 7]{Mordell} for finding parametric rational solutions of a homogeneous polynomial of degree $2$ in several variables.
\begin{Remark}
\label{Rem2}
The points $(t^2,d),((t+1)^2,e),((t+2)^2,f), ((t+3)^2,g)$, where $d,e,f,g\in \Q(t,p,q,w)$ are given as above, are rational points on the elliptic surface $y^2=ax^3+bx+c$, where $a,b,c$ are defined in (\ref{eqabc}).
\end{Remark}

Now, we assume that $((t+4)^2,h)$ is a rational point on the elliptic curve $y^2=ax^3+bx+c$. In particular, there exists a $5$-term sequence of consecutive squares on the latter curve. Then one has
\begin{eqnarray}
\label{eq:ellipticsurface}
h^2=A p^4+B p^3 +C p^2 +D p + E
\end{eqnarray}
with
{\footnotesize\begin{eqnarray*}
A&=&(140 + 246 t + 166 t^2 + 51 t^3 + 6 t^4)^2\\
  B&=&\frac{4(5 + 2 t)^2 (196 + 322 t + 202 t^2 + 57 t^3 +6 t^4) ( (87 + 83 t + 27 t^2 + 3 t^3)q -3 (52 + 58 t + 22 t^2 + 3 t^3) w)}{3 + 2 t}\\
  C&=&\frac{2(5 +2 t) }{3 + 2 t}(5  (78330 + 250402 t + 346118 t^2 + 271991 t^3 +132943 t^4 + 41217 t^5 + 7851 t^6 + 828 t^7\\ &+& 36 t^8)q^2-12  (41580 + 154358 t + 243358 t^2 + 217563 t^3 + 121708 t^4 + 43727 t^5 + 9852 t^6 + 1272 t^7 \\&+& 72 t^8)q w +15 (2 + t)^2 (-2828 - 2552 t + 784 t^2 + 2364 t^3 + 1395 t^4 +360 t^5 + 36 t^6) w^2)\\
  D&=&\frac{12 (35 + 24 t +4 t^2)}{3 + 2 t} ( (5655 + 17662 t + 23115 t^2 + 16782 t^3 + 7345 t^4 +1938 t^5 + 285 t^6 + 18 t^7)q^3\\ &+& (6660 + 18502 t + 22620 t^2 + 15567 t^3 + 6515 t^4 +1683 t^5 + 255 t^6 + 18 t^7)q^2 w -5  (3708 + 11842 t\\ &+& 16104 t^2 + 12237 t^3 + 5651 t^4 +1593 t^5 + 255 t^6 + 18 t^7)q w^2 + 3 (2 + t)^2 (260 + 888 t + 972 t^2 + 544 t^3\\ &+& 153 t^4 + 18 t^5) w^3)\\
        E&=&9 ( (65 + 141 t + 111 t^2 + 41 t^3 + 6 t^4)^2 q^4 +8  (22750 + 79965 t + 121251 t^2 + 105282 t^3 + 57708 t^4\\ &+&20529 t^5 + 4643 t^6 + 612 t^7 + 36 t^8) w q^3 -2 (120300 + 457050 t + 737244 t^2 + 678163 t^3 + 394077 t^4 \\&+&149001 t^5 + 35957 t^6 + 5088 t^7 + 324 t^8) w^2 q^2  +8  (2 + t)^2 (1750 + 6380 t + 7959 t^2 + 5294 t^3 + 1987 t^4 \\&+&408 t^5 + 36 t^6)q w^3 + (20 + 66 t + 66 t^2 + 31 t^3 +6 t^4)^2 w^4).
\end{eqnarray*}}

\begin{Theorem}
\label{thm1}
The curve $\CC:Y^2=A X^4+B X^3 +C X^2 +D X + E$ defined over $\Q(t)$ is birationally equivalent over $\Q(t,p,q,w)$ to an elliptic curve $\EE$  with $\rank \EE(\Q(t,p,q,w))\ge 1$.
\end{Theorem}
\begin{Proof}
After homogenizing the equation describing $\CC$, one obtains $Y^2=A X^4+B X^3Z +C X^2Z^2 +D XZ^3 + EZ^4$ with a rational point $R=(X:Y:Z)=(1:140 + 246 t + 166 t^2 + 51 t^3 + 6 t^4:0)$. The curve $\CC$ is birationally equivalent to the cubic curve $\EE$ defined by the equation $V^2=U^3-27 I U-27J$, \cite{StollCremona}, where $I=12AE-3BD+C^2$ and $J=72 ACE + 9 BCD -27 AD^2 -27B^2E -2 C^3$. The discriminant $\Delta(\EE)$ of $\EE$ is given by $(4I^3-J^2)/27$, and the specialization of $\EE$ is singular only if $\Delta(\EE)=0$.
Moreover, the point $P=\displaystyle \left( 3\frac{3B^2 -8AC}{4A}, 27\frac{B^3 + 8A^2D - 4ABC}{8A^{3/2}}\right)$ lies in $\EE(\Q(t,p,q,w))$ since $A$ is a square. One considers the specialization $\displaystyle t=1 , q=\frac{81}{40} ,w=1 $ to obtain the specialization $\displaystyle \widetilde{P}=\left(\frac{-4786935489}{100},\frac{-56568093052527}{50}\right)$ of the point $P$ on the specialized elliptic curve $$\widetilde{\EE}: y^2=x^3-\frac{147183268996968521373}{10000}x+\frac{171278570868444028577352480093}{250000}.$$ Using $\Magma$, \cite{Bosma}, the point $\widetilde{P}$ is a point of infinite order on $\widetilde{\EE}$. Therefore, according to Silverman's specialization Theorem, the point $P$ is of infinite order on $\EE$.
\end{Proof}

\begin{Corollary}
\label{cor1}
For any nontrivial sequence of consecutive rational squares $t_0^2,(t_0+1)^2,(t_0+2)^2,(t_0+3)^2, (t_0+4)^2$, there exist infinitely many elliptic curves $E_m:y^2=a_mx^3+b_mx+c_m,\;m\in\Z\setminus\{0\},$ such that $(t_0+i)^2,i=0,1,2,3,4,$ is the $x$-coordinate of a rational point on $E_m$. Moreover, these five rational points are independent.
\end{Corollary}
\begin{Proof}
We fix $t=t_0$, $q=q_0$, and $w=w_0$ in $\Q$. Substituting these values into (\ref{eq:ellipticsurface}), one obtains the elliptic curve
 \[\CC_{t_0,q_0,w_0}:h^2=Ap^4+Bp^3+Cp^2+Dp+E,\; A,B,C,D\in\Q,\] with positive rank, see Theorem \ref{thm1}. Now, one fixes a point $P=(p,h)$ of infinite order in $\CC_{t_0,q_0,w_0}(\Q)$. For any nonzero integer $m$, we set $mP=(p_m,h_m)$ to be the $m$-th multiple of the point $P$ in $\CC_{t_0,q_0,w_0}(\Q)$.

  Now, one substitutes $t=t_0,q=q_0,w=w_0$, and $p=p_m$ into the formulas for $d,e,f,g\in\Q(t,p,q,w)$ in (\ref{eq2}) in order to obtain the rational numbers $d_m,e_m,f_m,g_m$, respectively. Then one substitutes $d_m,e_m,f_m$ into the formulas for $a,b,c\in\Q(t,d,e,f)$ in (\ref{eqabc}) to get the rational numbers $a_m,b_m,c_m$, respectively.

  To sum up, one constructed an infinite family of elliptic curves $E_m:y^2=a_mx^3+b_mx+c_m$, where $m$ is a nonnegative integer. The latter infinite family $E_m$ of elliptic curves satisfies the property that the points $(t_0^2,d_m),((t_0+1)^2,e_m),((t_0+2)^2,f_m),((t_0+3)^2,g_m),((t_0+4)^2,h_m)\in E_m(\Q)$. Thus, one obtains an infinite family of elliptic curves with a $5$-term sequence of rational points whose $x$-coordinates form a sequence of consecutive squares in $\Q$.

  To show that the points $(t_0^2,d_m),((t_0+1)^2,e_m),((t_0+2)^2,f_m),((t_0+3)^2,g_m),((t_0+4)^2,h_m)\in E_m(\Q)$ are independent, one specializes $\displaystyle t=1 , q=81/40 ,w=1 $ which yields the existence of the infinite point $\displaystyle (p,h)=\left(\frac{2201}{2320},\frac{-62736289}{18852320}\right)\in \CC_{1,81/40,1}(\Q)$. Therefore, the specialization $t=1,q=81/40,w=1,p=2201/2320$ gives us the specialized elliptic curve
 \[E_1: {\footnotesize y^{2}=\frac{42674183}{52786496000}x^3-\frac{612989889 }{7540928000}x+\frac{1180698375893607 }{2487869785676800}}\]
                                 with the following set of rational points in $E_1(\Q)$:
                                 \[\left(1,\frac{-2367005}{3770464}\right), \left(2^2,\frac{8455597}{18852320}\right), \left(3^2,\frac{-10868031}{18852320}\right),\\ \left(4^2,\frac{-29720351}{18852320}\right), \left(5^2,\frac{-62736289}{18852320}\right).\] Using $\Magma$, \cite{Bosma}, these rational points are independent.

  According to Silverman's Specialization Theorem, it follows that the points $(t_0^2,d_m),((t_0+1)^2,e_m),((t_0+2)^2,f_m),((t_0+3)^2,g_m),((t_0+4)^2,h_m)$ are independent in $E_m$ over $\Q(t_0,q_0,w_0,p_m)$.
\end{Proof}
\begin{Remark}
Corollary \ref{cor1} implies the existence of an infinite family of elliptic curves whose rank $r\ge 5$.
\end{Remark}

\begin{Remark}
One notices that a sequence of consecutive squares on an elliptic curve gives rise to a set of rational points on some hyperelliptic curve of genus $2$, see the proof of Proposition \ref{prop1}. Therefore, according to Corollary \ref{cor1}, we are able to construct an infinite family of hyperelliptic curves $\CC$ such that $|\CC(\Q)|\ge 5$.
\end{Remark}
\hskip-12pt\emph{\bf{Acknowledgements.}}
We would like to thank Professor Nabil Youssef, Cairo University, for his support, thorough reading of the manuscript, and several useful suggestions.

\end{document}